\documentclass[aip,cha,graphicx,reprint]{revtex4-1}

\usepackage{amsmath}
\usepackage{amssymb}
\usepackage{graphicx}
\usepackage{dcolumn}
\usepackage{bm}

\usepackage[utf8]{inputenc}
\usepackage[T1]{fontenc}
\usepackage{mathptmx}
\usepackage{etoolbox}
\usepackage[colorlinks=true]{hyperref}
\usepackage{xcolor}
\usepackage[english]{babel}
\usepackage{subfig}
\usepackage{booktabs}
\usepackage{braket}

\newcommand{\E}{\mathbb{E}} 
\newcommand{\ergset}{\mathcal{X}} 
\DeclareMathOperator*{\argmin}{arg\,min}

\makeatletter
\def\@email#1#2{%
 \endgroup
 \patchcmd{\titleblock@produce}
  {\frontmatter@RRAPformat}
  {\frontmatter@RRAPformat{\produce@RRAP{*#1\href{mailto:#2}{#2}}}\frontmatter@RRAPformat}
  {}{}
}%
\makeatother
\begin{document}

\title{Computing Chaotic Time-Averages from Few Periodic or Non-Periodic Orbits}
\author{Joshua L. Pughe-Sanford}
\email{jpughesanford@gatech.edu}
\author{Sam Quinn}
\author{Teodor Balabanski}
\author{Roman O. Grigoriev}
\affiliation{School of Physics, Georgia Institute of Technology, 837 State St NW, Atlanta, GA 30332, USA}

\date{\today}

\begin{abstract}
For appropriately chosen weights, temporal averages in chaotic systems can be approximated as a weighted sum of averages over reference states, such as unstable periodic orbits. Under strict assumptions, such as completeness of the orbit library, these weights can be formally derived using periodic orbit theory. When these assumptions are violated, weights can be obtained empirically using a Markov partition of the chaotic set. Here, we describe an alternative, data-driven approach to computing weights that allows for an accurate approximation of temporal averages from a variety of reference states, including both periodic orbits and non-periodic trajectory segments embedded within the chaotic set. For a broad class of observables, we demonstrate that weights computed with the proposed method significantly outperform those derived from periodic orbit theory or Markov models, achieving superior accuracy while requiring far fewer states---two critical properties for applications to high-dimensional chaotic systems.
\end{abstract}

\maketitle

\begin{quotation}
Chaotic systems, which arise in a myriad of natural and engineered settings, challenge traditional methods of predicting long-term behavior due to their sensitivity to initial conditions. While periodic orbit theory provides a foundational framework for approximating time-averaged observables in chaotic systems, its reliance on complete sets of unstable periodic orbits (UPOs) often limits its applicability in high-dimensional settings where such complete sets may be very expensive or even impossible to identify. This work introduces a data-driven, interpretable approach leveraging reproducing kernel Hilbert space (RKHS) interpolation, termed least-squares weighting (LSW), to estimate chaotic averages using arbitrary solution segments, including non-periodic "snippets." The LSW method achieves superior accuracy and faster convergence compared to state-of-the-art approaches such as periodic orbit theory, even with a limited library of states. By expanding the toolkit for analyzing high-dimensional chaos, our results open pathways for applications in fluid dynamics, magnetohydrodynamics, optics, and beyond, offering robust predictions in systems where traditional methods falter. 
\end{quotation}

\section{Introduction}

Deterministic chaos is pervasive in natural and engineered systems, with turbulence in fluids \cite{christensen2021}, plasma \cite{yamada2008} and nonlinear optics \cite{picozzi2014} representing some well-known and practically important examples of high-dimensional chaotic systems. Despite being governed by deterministic equations, chaotic systems 
can only be described dynamically 
on short temporal intervals  due to their inherent sensitivity to initial conditions \cite{lorenz1963}.  Hence, on long time scales, one is forced to resort to a statistical description, with a key objective to compute temporal averages of various physical observables. 
An illustrative example of dynamical vs. statistical descriptions of chaotic systems is the prediction of weather vs. the prediction of climate.

Here we aim to predict the averages of
physical observables, $a({\bf x})$, defined as functions that map the instantaneous state of the system, ${\bf x}\in\mathcal{M}$, to a real number. So long as the dynamics are ergodic---a property many chaotic systems either have or may be presumed to have in practice---Birkhoff's ergodic theorem  \cite{birkhoff1931} ensures that temporal averages may be computed
as the expectation value 
over a Sinai-Ruelle-Bowen (SRB) measure $\mu$ \cite{Sinai1972, Ruelle1976, Bowen2008},
\begin{align}
     \E[a]  &= \lim_{T\to\infty}\frac{1}{T}\int_0^{T} {a}(\mathbf{x}(t))\,dt 
     \label{eq:tave}\\
     &\equiv \int_\ergset  {a}(\mathbf{x})\,d\mu(\mathbf{x}) \label{eq:mave},
\end{align}
for all ${a}({\bf x})$ that are $\mu$-integrable. 
Here, ${\ergset}$
denotes an ergodic subset of the state space $\mathcal{M}$, and ${\bf x}(t)$ is a trajectory confined to, or approaching, ${\ergset}$.
Note that $\E[a]$ is 
invariant with respect to the choice of ${\bf x}(0)\in \ergset$.

Evaluating $\E[\,\cdot\,]$ exactly requires an infinitely long trajectory ${\bf x}(t)$ that explores $\ergset$ 
or the knowledge of $\mu$ itself, either of which is only available in highly exceptional cases \cite{grigoriev1998}.
The most direct and common approach is to approximate $\E[\,\cdot\,]$ using a finite but sufficiently long trajectory. However, in high-dimensional systems, storing such trajectories is exceedingly expensive and cumbersome: for instance, a turbulent fluid flow on a relatively small computational domain and over a relatively short time interval requires hundreds of terabytes of storage \cite{graham2016}.

A more practical, and far more data-efficient, alternative is to define a collection $\{ \mu_p \}_{p=1}^P$ of measures supported on subsets $S_p \equiv \operatorname{supp} \mu_p \subseteq \mathcal{M}$, such that
\begin{align}\label{eq:genexpansion}
\hat{\E}[a] \equiv \sum_{p=1}^P w_{p} \, \E_p[a]  \approxeq \E[a],
\end{align}
is a good approximation of $\E[a]$, where 
\begin{align}
\E_p[a] = \int_{S_p} {a}(\mathbf{x})d\mu_{p}({\bf x})
\end{align}
is the integral of $a$ over $\mu_p$, and $w_{p}$ are scalar weights.
Various expansions of this form 
have been proposed \cite{ulam2004, chaosbook, fernex2021}, with different choices of $S_p$, means for computing $w_{p}$, advantages and drawbacks. 

\autoref{eq:genexpansion} can only be used to predict the averages of observables for which $\E[a]$ and $\E_p[a]$ exist and are finite. It is common to require that all measures are normalizable, such that $0<|\E[1]|, |\E_p[1]| <\infty$. This condition ensures that $\E[a]$ and $\E_p[a]$ are finite for all bounded observables, i.e., those satisfying $\sup_{\bf x\in\mathcal{M}}|a({\bf x})|<\infty$. Without loss of generality, it is further assumed that all averages are normalized, $\E[1] = \E_p[1] =1.$

An expansion such as \autoref{eq:genexpansion} is particularly well suited for computing the averages that depend on hyper-parameters. Examples include the $\gamma$-norm of the state vector, $a_\gamma({\bf x}) = \left(\sum_i x_i^\gamma\right)^{-\gamma}$, the variance of the flow about a reference state, $a_{\boldsymbol{\gamma}}({\bf x}) = ({\bf x}-\boldsymbol{\gamma})^\top ({\bf x}-\boldsymbol{\gamma})$, or even the Fourier Transform of the SRB measure, $a_{\boldsymbol{\gamma}}({\bf x}) = e^{i\boldsymbol{\gamma}\cdot {\bf x}}$. The nonlinear dependence of $a_{\boldsymbol{\gamma}}$ on $\boldsymbol{\gamma}$ renders direct computation of 
the averages of such an observable over a range of $\gamma$
very expensive. Traditionally, the hyper-parameter space must be discretized, and an infinite-time average $\E[a_{\boldsymbol{\gamma}}]$ computed for each $\boldsymbol{\gamma}$ on a grid. \autoref{eq:genexpansion} allows for, at every gridpoint, the expensive average over the chaotic trajectory to be replaced with inexpensive averages over reference states. For fine discretizations of the $\boldsymbol{\gamma}$, this can result in massive speedup. 

A particularly desirable choice for the subsets $S_p$ is the unstable periodic orbits (UPOs) \cite{eckhardt1994, cvitanovic1995,cvitanovic1999}. Periodic orbits are special, repeating solutions of a dynamical system satisfying ${\bf x}_p(t) = {\bf x}_p(t+T_p)$, where $T_p$ is the period of motion. Due to periodicity, averages over UPOs can be computed exactly in finite time,
\begin{align}\label{eq:poave}
    \E_p[a] &= \frac{1}{T_p}\int_0^{T_p} {a}(\mathbf{x}_p(t))\ dt. 
\end{align} 
allowing the infinite time-average present in \autoref{eq:genexpansion} to be replaced with a collection of finite time-averages. 

This choice, which is the basis of the Periodic Orbit Theory (POT), leverages the fact that periodic orbits are dense in the closure of the ergodic set of Axiom A systems \cite{eckmann1985, smale1967}. In particular, for Axiom A systems, it can be proven that \autoref{eq:genexpansion} converges to equality for smooth observables \cite{chaosbook}  in the limit $P\to \infty$. Moreover, the weights $w_p$ may be computed analytically in terms of the stability properties of the 
UPOs \cite{chaosbook_23.5}. A straightforward generalization of POT exists for systems with continuous symmetries, where UPOs become relative \cite{budanur2015}; we will refer to both types of solutions as simply UPOs.

Periodic orbits often exemplify the same spatiotemporal features as the chaotic flow itself \cite{waleffe2001,toh2003,chantry2014,shekar2018} and are considered to form a dynamical skeleton of $\mathcal{X}$ \cite{budanur2017,crowley2023}. Hence, when periodic orbits are known, they generate a remarkably computationally efficient and \textit{interpretable} representation of the chaotic averages, particularly in high dimensions. For this reason, practitioners have sought, and struggled, to apply POT to high-dimensional chaotic flows for more than two decades \cite{kawahara2001,chandler2013,suri2020,page2024,yalniz2021,chaosbook}. 

It is notoriously difficult to compute a sufficiently comprehensive collection of periodic orbits in high dimensions \cite{chandler2013}. For instance, in turbulent fluid flows, the number of known periodic orbits rarely exceeds a few tens \cite{budanur2017,krygier2021,crowley2023}. The orbits that are discovered tend to be poorly distributed within the chaotic set, a limitation that only state-of-the-art methods have begun to address \cite{page2024}. 

Also, the convergence of POT relies critically on the existence of a \textit{symbolic dynamics}. Symbolic dynamics is a mathematical framework that simplifies the description of a system by partitioning its state space into distinct regions and representing trajectories by the sequence of regions they visit \cite{lind1995}. In ergodic systems, this partition is meticulously constructed to reflect the underlying topology of the chaotic set. 
For periodic orbits, which repeat in time, their symbolic sequences are likewise periodic, with the symbol length defined as the number of symbols before the sequence repeats. A library of orbits is said to be complete up to symbol length $L$ if it includes---and only includes---all periodic orbits with symbol length $l<L$. POT is proven to converge exponentially fast for successively \textit{complete} libraries, in the limit $L\to\infty$. However, POT is only as accurate as the shortest orbit truncated from the expansion; if a short orbit is missed, adding more orbits only negligibly improves the accuracy of the method \cite{chaosbook_23.1}. 

It is important to note that most systems of practical importance do not admit or do not have known symbolic dynamics, rendering it impossible to tell when a library is complete and further complicating the usage of POT. Stability ordering---an extension of POT to incomplete libraries and systems without symbolic dynamics---has been proposed \cite[Ch. 23.7]{chaosbook}. However, this is an active area of research and the expected accuracy of stability ordering has not yet been established. These issues, collectively, have prevented POT from producing accurate predictions in high-dimensional systems \cite{chandler2013}. 

A popular alternative to POT 
is Ulam's method \cite{ulam2004, maiocchi2022}, which typically partitions a compact set enclosing $\mathcal{X}$ into disjoint subsets $\{ S_p \}_{p=1}^P$, each imbued with a uniform measure. This approach is proven to approximate averages over any absolutely continuous measure in the limit that the size of the subsets becomes vanishingly small. However, this method also struggles in applications to high-dimensional flows. Even for a moderate number of dimensions, Ulam's method quickly becomes intractable due to the proliferation of sets and difficulty associated with integrating over the uniform measure in $d$ dimensions. Some studies are attempting to address these issues with novel clustering methods \cite{fernex2021}.

Lacking many prerequisites for a working periodic orbit theory,
practitioners working in high-dimensional systems have pursued a periodic orbit decomposition of chaotic averages in which---in lieu of the POT weights---$w_p$ are computed from a Markov process defined on the Voronoi volumes about each UPO \cite{yalniz2021, page2024}. Empirically, Markov-based approaches achieve moderate accuracy on the order of a few percent. 
However, the Markov-based approach lacks theoretical guarantees regarding the types of observables that can be accurately predicted using small sets of UPOs and whether it converges in the limit $P\to\infty$. 

In this letter, we introduce an interpretable data-driven approach to determining optimal weights for a finite, and possibly quite small, library of measures $\{ \mu_p \}_{p=1}^P$.
It utilizes Kriging \cite{Cressie1990}, a method that leverages a Reproducing Kernel Hilbert Space (RKHS) to perform interpolation and regression on abstract data. Kriging is known for its accuracy in applications to high-dimensional data sets \cite{Wang2024, Nan2024, Song2024}. 
The method enables averages over a target measure $\mu$ to be estimated from \textit{any} choice of reference measures $\{ \mu_p \}_{p=1}^P$. It also provides rigorous claims as to which observables averages can be computed with zero error. The proposed method finds multiple distributions of weights that are more accurate and converge more rapidly in $P$ than both POT and Markov models. 

\section{The Method} 

The method being proposed here is quite simple. To illustrate it, we begin with a tutorial before presenting a more rigorous derivation.

\subsection{Tutorial\label{sec:tutorial}}

Suppose that $P=2$ and that $S_p=\{{\bf x}_p(t)\ |\ 0\le t<T_p\}$ where $p=1,2$ are sets associated with periodic orbits. Then, \autoref{eq:genexpansion} reduces to
\begin{equation}\label{eq:tutorial_discrete}
    \E[a] \approxeq w_1 \E_1[a]+w_2 \E_2[a].
\end{equation}
The above equation is linear, with two unknowns $w_1$ and $w_2$. By evaluating \autoref{eq:tutorial_discrete} on two or more test observables, we may solve for the weights using the method of least squares. 

For instance, suppose that $a_1({\bf x})$ and $a_2({\bf x})$ are observables of particular interest. Then, weights that best predict the chaotic averages of $a_1$ and $a_1$ from averages over orbits $1$ and $2$ will approximately satisfy \autoref{eq:tutorial_discrete} for both $a_1$ and $a_2$:
\begin{equation}\label{eq:linear_system}
\begin{pmatrix} \E[a_1] \\ \E[a_2] \end{pmatrix}\approx \begin{pmatrix} \E_1[a_1] & \E_2[a_2] \\ \E_1[a_2] & \E_2[a_2]\end{pmatrix}\begin{pmatrix} w_1 \\ w_2 \end{pmatrix}.
\end{equation}
A least squares solution 
\begin{equation}\label{eq:solution_discrete}
\begin{pmatrix} w_1 \\ w_2 \end{pmatrix} = \begin{pmatrix} \E_1[a_1] & \E_2[a_1] \\ \E_1[a_2] & \E_2[a_2]\end{pmatrix}^\dagger\begin{pmatrix} \E[a_1] \\ \E[a_2] \end{pmatrix}
\end{equation}
always exists, where $[\ \cdot \ ]^\dagger$ indicates a pseudo-inverse. Solving \autoref{eq:linear_system} when the matrix in non-invertible is further discussed in \autoref{sec:uniquenessandconstraints}.

\begin{figure*}
    \centering
    \subfloat[]{\includegraphics[width=0.292\linewidth]{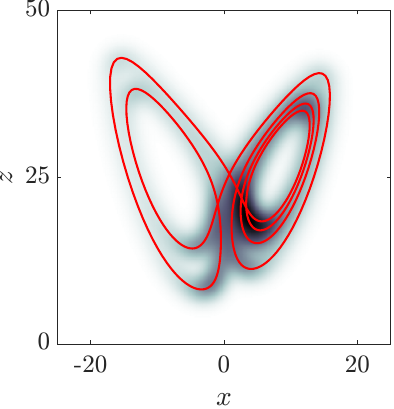}}\hspace{5mm}
    \subfloat[]{\includegraphics[width=0.25\linewidth]{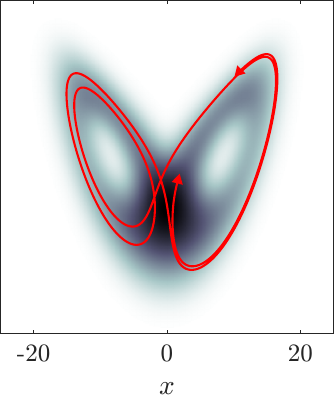}}\hspace{5mm}
    \subfloat[]{\includegraphics[width=0.25\linewidth]{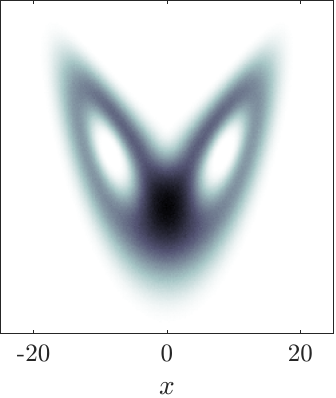}}
    \caption{An illustration of (a) periodic orbit density $f_{p}({\bf x})$ and (b) chaotic orbit density $f({\bf x})$ over the Lorenz 1963 dynamics using a Gaussian kernel with unit variance. In each, the trajectory that generates the data is shown in red. Darker colors indicate larger values of $f$. Red arrows illustrate that the chaotic trajectory extends infinitely in both directions, without repeating.  (c) A histogram of all states visited by a long chaotic trajectory, color indicated the probability that a bin is visited by the chaotic trajectory. Darker bins are visited more often. In all panels, the function (or histogram) is averaged over $y$ for plotting purposes. }
    \label{fig:densities}
\end{figure*}

As written, the weights given in \autoref{eq:solution_discrete} are heavily biased towards predicting well the averages of $a_1$ and $a_2$. They may predict other observables' averages arbitrarily poorly. To avoid biasing the least squares solution towards any particular physical observable, consider instead a spatially localized observable, such as a Gaussian kernel
\begin{equation}
G_\theta({\bf x},{\bf y}) = \exp\left[-\frac{|{\bf y}-{\bf x}|^2}{2\theta}\right],
\end{equation}
with variance $\theta$ and mean ${\bf y}$. Evaluating \autoref{eq:tutorial_discrete} on $G_\theta$, it follows that 
\begin{equation}\label{eq:tutorial_continuous}
f({\bf y}) \approxeq w_1f_1({\bf y}) + w_2 f_2({\bf y}),
\end{equation}
where 
\begin{align}
    f({\bf y}) &= \E[G_\theta(\ \cdot \ ,{\bf y})]\\
    &= \int_\mathcal{M} G_\theta({\bf x},{\bf y}) d\mu({\bf x})\\
    &= \lim_{T\to\infty}\frac{1}{T}\int_0^TG_\theta({\bf x}(t),{\bf y}) dt 
\end{align}
is a function related to the probability that the chaotic trajectory visits the $\theta$-neighborhood of ${\bf y}$, and \begin{align}
    f_p({\bf y}) &= \E_p[G_\theta(\ \cdot \ ,{\bf y})]\\
    &= \int_{S_p} G_\theta({\bf x},{\bf y}) d\mu_p({\bf x})\\
    &= \frac{1}{T_p}\int_0^{T_p}G_\theta({\bf x}_p(t),{\bf y}) dt 
\end{align}
is a function related to the probability that orbit $p$ visits the $\theta$-neighborhood of ${\bf y}$. In the limit that $\theta$ becomes small, $f({\bf x})$ limits to---up to normalization---the probability that a randomly sampled chaotic state is ${\bf x}$. Indeed, $\theta$ plays a regularizing role quite similar to that of noise; it has been shown previously that---in the presence of noise---this probability can be expanded as a linear sum of Gaussians supported on periodic points in the chaotic set \cite{heninger2015}. \autoref{eq:tutorial_continuous} expresses a similar idea.  Examples of $f({\bf x})$ and $f_p({\bf x})$ are illustrated in \autoref{fig:densities} for the Lorenz 1963 system \cite{lorenz1963} with canonical parameter values $\sigma=10$, $\rho=28$, and $\beta=8/3$; Lorenz is a reduced order model for atmospheric convection 
and a prototypical example of deterministic chaos.

Since all choices of ${\bf y}\in \mathcal{M}$ constrain the weights, \autoref{eq:tutorial_continuous} represents an \textit{over-constrained} system of equations for weights $w_1$ and $w_2$. It can also be solved using the method of least squares. The resulting optimal weights are:
\begin{equation}\label{eq:solution_continous}\begin{pmatrix} w_1 \\ w_2 \end{pmatrix} = \begin{pmatrix} \langle f_1, f_1\rangle  & \langle f_1, f_2\rangle \\ \langle f_2, f_1\rangle & \langle f_2, f_2 \rangle \end{pmatrix}^\dagger\begin{pmatrix} \langle f_1, f \rangle\\ \langle f_2, f\rangle \end{pmatrix},
\end{equation}
where 
\begin{equation}\label{eq:inner_product}
    \langle f , g\rangle = \int_\mathcal{M} f({\bf x})g({\bf x})d{\bf x}
\end{equation}
is an inner product on $\mathcal{M}$. 

Integration over a state space---as required by \autoref{eq:inner_product}---is likely numerically intractable. However, a well chosen kernel will allow the integral over $\mathcal{M}$ to be taken analytically. Notice that for a Guassian kernel in $\mathcal{M}=\mathbb{R}^d$, only integrals over orbits $p$ and $q$ are required:
\begin{align}
    \langle f_p, f_q\rangle &= \int_\mathcal{M} f_p({\bf x}) f_q({\bf x})\ d{\bf x} \\
    &= \int_{S_p}d\mu_p({\bf y})\ \int_{S_q}d\mu_q({\bf z})\ \int_\mathcal{M} G_\theta({\bf y},{\bf x})G_\theta({\bf x},{\bf z})\ d{\bf x}\nonumber \\
    &= (\pi\theta)^{d/2} \int_{S_p}d\mu_p({\bf y})\ \int_{S_q}d\mu_q({\bf z})\ G_{2\theta}({\bf y},{\bf z})\label{eq:tutorialip}\\
    &=\frac{(\pi\theta)^{d/2}}{T_pT_q} \int_0^{T_p}\int_0^{T_q} G_{2\theta}({\bf x}_p(t),{\bf x}_q(\tau)) dt d\tau.\label{eq:overorbit}
\end{align}
Hence, these inner products can be computed as averages over orbits, which are computationally inexpensive.

Unlike the weights given in \autoref{eq:solution_discrete}, the solution provided in \autoref{eq:solution_continous} is unbiased towards any physical observable. Rather, these weights optimize the global reconstruction of $f({\bf x})$ from functions $f_p({\bf x})$. 

In the next section, we re-derive the method for an arbitrary kernel function observable and define the subspace of observables $\mathcal{A}$ whose averages are predicted by the method with zero error.

\subsection{\label{sec:derivation}General Approach}

Consider a smooth, bounded kernel function $k({\bf x},{\bf y}) \in L^2(\mathcal{M}\times\mathcal{M})$ satisfying 
\begin{align}\label{eq:G2}
    \int_{\mathcal{M}} k({\bf x},{\bf z})k({\bf z},{\bf y})\ d{\bf z} = k^2({\bf x},{\bf y}),
\end{align}
where $k^2$ is an arbitrary induced kernel function.
Because $k$ is smooth and bounded, so is $k^2$, and both are guaranteed to be $\mu$- and $\mu_p$-integrable. This allows us to define the correlation between arbitrary normalizable measures $\mu'$ and $\mu''$ (supported on sets $S', S''\subset \mathcal{M}$) as
\begin{align}\label{eq:measkern}
    K(\mu',\mu'') &\equiv \int_{\mathcal{M}} \mathbb{E}_{\mu'}\left[k({\,\cdot\,},{\bf z})\right]\mathbb{E}_{\mu''}\left[k({\,\cdot\,},{\bf z})\right]\ d{\bf z}\nonumber\\
    &= \int_{S'} d\mu'({\bf x}) \int_{S''} d\mu''({\bf y})\ k^2({\bf x},{\bf y}).
\end{align}
Intuitively, one can think of $K({\mu'},{\mu''})$ as an inner product of measures $\mu'$ and $\mu''$. Indeed, when $k$ is a Guassian kernel function,
$K(\mu_p,\mu_q) = \langle f_p,f_q\rangle$
is the inner product from \autoref{eq:tutorialip}.

Let us associate a unique observable $a_p({\bf x})$ with each orbit measure $\mu_p$ by defining
\begin{equation}\label{eq:defnap}
    a_p({\bf x})\equiv \int_{S_p} d\mu_{p}({\bf y})\ k^2({\bf y},{\bf x}).
\end{equation}
Let $\mathcal{A} = \textrm{span}\left(a_1({\bf x}),\cdots,a_P({\bf x})\right)$ be the span of these functions. For any observable $a\in\mathcal{A}$, i.e.,
\begin{equation}\label{eq:ainA}
a({\bf x}) = \sum_{p=1}^P c_p a_p({\bf x}),
\end{equation}
using \autoref{eq:measkern} we find that
\begin{equation}\label{eq:constraints}
\int d\mu'({\bf x})\ a({\bf x}) = \sum_p c_p K(\mu_p,\mu').
\end{equation}
Evaluating \autoref{eq:constraints} at $\mu' = \mu_q$ for $q = 1,\dots,P$, we find:
\begin{align}
\sum_{p=1}^P c_p A_{pq} = \E_{q}[a],
\end{align}
where $A_{pq} = K(\mu_p,\mu_q)$. 

Define weights $w_p$ as the solution to the system 
\begin{align}\label{eq:rkhsw}
\sum_{q=1}^P A_{pq}w_q = \E[a_p].
\end{align}
Then, averaging both sides of \autoref{eq:ainA} over the SRB measure, it follows that
\begin{align}\label{eq:rkhsa}
\E[a] &= \sum_{p=1}^P c_p \E[a_p] \nonumber\\
&=\sum_{q=1}^P w_q \E_q[a]
=\hat{\E}[a],
\end{align}
where the equality between $\E[a]$ and $\hat{\E}[a]$ is {\it exact}, rather than approximate, for all ${a}\in \mathcal{A}$. 

Hence, \autoref{eq:rkhsw} recasts computing optimal weights as solving a simple linear system of the form $A\vec{w} = \vec{b}$, where $b_p = \E[a_p]$ is a vector of $P$ time-averages. We call $\hat{\E}[a]$ with the weights given by \autoref{eq:rkhsw} a least-squares weighted (LSW) estimate of $\E[a]$. 

For observables $a\notin\mathcal{A}$, using the Least Squares weights to compute $\hat{\E}[a]$ may not yield an exact value for $\E[a]$. 
However, the error can remain quite small for observables outside $\mathcal{A}$. 
For instance, since the observables in $\mathcal{A}$ are smooth functions of the state ${\bf x}$ (as a consequence of smoothness of the kernel), one can assume that low error can be similarly achieved for smooth observables. We indeed find this to be the case, as evidenced in \autoref{sec:results}.

\section{Results\label{sec:results}}

We compare the accuracy of different approaches using the Lorenz 1963 system. Lorenz admits a symbolic dynamics, making it an ideal testbed in which to compare LSW against POT. For an introduction to symbolic dynamics in Lorenz, we point the reader to \citet{viswanath2003}. In this system, an orbit's symbol length is equal to the number of times it winds around the left and right sides of chaotic set, before repeating; see \autoref{fig:densities}(a) for an example of a UPO with symbol length $l=6$. We will use the $P\leq 125$ UPOs of shortest symbol length from a collection of more than a thousand found by Viswanath \cite{viswanath2003, viswanathdata}. To ensure a fair and rigorous comparison between LSW and POT, we compute POT weights exclusively for complete libraries of symbol length $l = 2,\dots,9$. The sizes of these complete libraries are $P=1,3,6,12,21,39,69$ and $125$, respectively. As previously stated, we are particularly interested in computing averages when $P$ is small, which will be the most realistic regime in applications to higher dimensional chaotic systems. 

In addition to computing averages for periodic orbits, we will also consider a 
short, finite-time intervals sampled from chaotic trajectory (which we call ``snippets"). Computing periodic orbits requires expensive numerical solvers. In contrast, sampling snippets from numerical integration is trivial; even computing snippets that are \textit{almost}-periodic from auto-recurrence is comparatively inexpensive. Snippets allow us to explore the applicability of a LSW estimate to systems for which no, or too few, UPOs are available.

To generate snippets, we compute a chaotic trajectory that explores $\ergset$ for $\sum_{p=1}^{125} T_p$ time units and, for simplicity, divide it uniformly into $125$ snippets of equal duration. In this way, the total length of all the UPOs and all snippets considered is the same. We do not ensure that the snippets are well-distributed on the attractor or enforce upon the collection of snippets any symmetries of the governing equations.

\begin{figure}
    \center
    \includegraphics[width=0.35\textwidth]{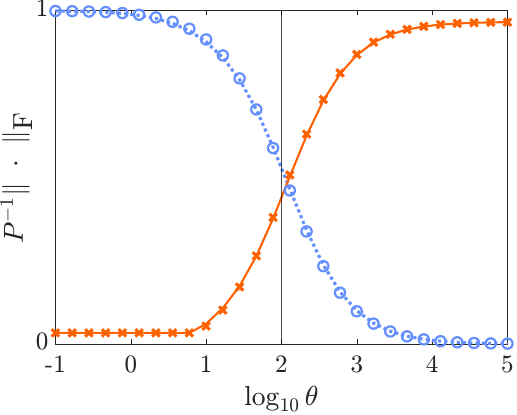}
    \caption{\label{fig:choiceoftheta}
    The Frobenius distance between $A$ and the matrix of all ones (blue, dotted line, o's) or the identity matrix (orange, solid line, x's) for $P=30$. Lines (symbols) correspond to matrices $A$ computed using UPOs (snippets). Black line indicates the value of $\theta=10^2$ used in \autoref{fig:scaling} and \autoref{tab:averages}. 
    }
\end{figure}

We use a Gaussian kernel function $k=G_\theta$, as introduced in the tutorial. For $\theta$ too large or too small, the correlation matrix $A$ reduces to a matrix of all ones or to a diagonal matrix,
respectively; in the former limit, the matrix $A$ becomes entirely singular and, in the latter, the LSW weights limit to Markovian-like  weights defined on sets $S_p$ (that is, $w_p \to \mathbb{E}\left[I({\bf x},S_p)\right]$, where $I({\bf x},S_p)$ is an indicator function equal to unity when ${\bf x}\in S_p$ and zero otherwise). Both extremes wash away the rich correlations between measures $\mu_p$. We use $\theta=10^2$ which, consistently, across both UPOs and snippets, is roughly equidistant---in terms of the Frobenius norm---from either extreme, as illustrated by \autoref{fig:choiceoftheta}.
For reference, the corresponding value of $\sqrt{\theta}$ is of order $\sqrt{2\beta(\rho-1)}=12$ which sets the size of the attractor.
There are examples of more complex choices of $\theta$ in the literature; for instance, \citet{lippolis2010} define a spatially varying Gaussian kernel width, proportional to the local response of a dynamics to noise. Here, we find a global choice of $\theta$ based on the Frobenius distance
to be simple yet efficient. 

For any fixed $\theta$, $A$ is likely to become increasingly singular as $P$ increases, due to increasing overlap between the functions $a_p$. To penalize marginal modes from appearing in the solution weights, we solve \autoref{eq:rkhsw} using Tikhonov regularization:
\begin{equation}\label{eq:sys_tikhonov}
    \left(A+\alpha I\right)\vec{w} = \vec{b},
\end{equation}
with $\alpha = 10^{-10}$. The marginal modes of $A$ are discussed further in \autoref{sec:uniquenessandconstraints}. 

We compute the elements of matrix $A$ and $\vec{b}$ using the formulae:
\begin{align}
    A_{pq} &= \int_0^{T_p} \frac{dt}{T_p} \int_0^{T_q} \frac{d\tau}{T_q} \exp\left[ -\frac{\|{\bf x}_p(t)-{\bf x}_q(\tau)\|^2}{4\theta}\right] \label{eq:Kpq} \\
    b_q &= \lim_{T\to\infty}\int_0^{T} \frac{dt}{T} \int_0^{T_q} \frac{d\tau}{T_q} \exp\left[ -\frac{\|{\bf x}(t)-{\bf x}_q(\tau)\|^2}{4\theta}\right],\label{eq:Eap}
\end{align}
The largest computational cost associated with this method is the evaluation of the averages in $\vec{b}$. However, these averages are independent and can be computed in parallel (across orbits and in time) quite quickly, once and for all, no matter how many different averages we might need to estimate later.

We cross-validate the accuracy of LSW on a set of observables not in $\mathcal{A}$, namely,
$\mathcal{B} = \left\{1,x,y,z,x^2,xy,xz,y^2,yz,z^2\right\}$. All second order polynomials, including observables like $\dot{\bf x}$, lie in the span of the elements of $\mathcal{B}$. We define the relative error
\begin{equation}\label{eq:erel}
    E_{\text{rel}}(a) = \frac{\left|\E[a] - \hat{\E}[a]\right|}{\sqrt{\text{var}(a)}},
\end{equation}
and maximal error
\begin{equation}\label{eq:emax}
    E_{\text{max}} = \max_{b\in\mathcal{B}} E_\text{rel}(b),
\end{equation}
where the variance, $\text{var}(a)$, sets a natural scale for deviations from $\E[a]$. For observables with zero variance, we set $\text{var}(a)\equiv 1$. 

We compare the accuracy of the proposed approach against three alternatives. The first one is POT, which has been applied to Lorenz previously and show ton provide accurate estimates of averages for a variety of observables \cite{eckhardt1994}. The second is a Markov weighting. Lastly, we consider a uniform weighting, where $w_p = 1/P$. Uniform weighting is the simplest option which has been used previously as a baseline with which to compare other weighting schemes \cite{chandler2013}. While POT can only be applied to UPOs, all other methods can be applied to snippets as well.

\begin{table}
    \centering
    \begin{tabular}{c|cccc|ccc}
    & & \multicolumn{2}{c}{Orbits} & & \multicolumn{3}{c}{Snippets} \\ 
 & \textcolor[HTML]{DC267F}{POT} & \textcolor[HTML]{FFB000}{Uniform} & \textcolor[HTML]{785EF0}{Markov} & \textcolor[HTML]{648FFF}{LSW} & \textcolor[HTML]{FFB000}{Uniform} & \textcolor[HTML]{785EF0}{Markov} & \textcolor[HTML]{648FFF}{LSW} \\
 \hline
 \hline
 1 & $-\boldsymbol{\infty}$ & $-\boldsymbol{\infty}$ & $-\boldsymbol{\infty}$ & $-3.8$ & $-\boldsymbol{\infty}$ & $-\boldsymbol{\infty}$ & $-4.1$ \\
  $x$   & $-\boldsymbol{\infty}$ & $-1.4$ & $-1.7$ & $-3.2$ & $-1.2$ & $-1.6$ & $-3.1$ \\
 $y$    & $-\boldsymbol{\infty}$ & $-1.4$ & $-1.7$ & $-3.3$ & $-1.3$ & $-1.6$ & $-3.1$ \\
 $z$    & $-1.7$ & $-1.8$ & $-2.0$ & $-3.5$ & $-2.1$ & $-2.2$ & $-3.2$ \\
 $x^2$  & $-2.2$ & $-2.2$ & $-2.4$ & $-3.2$ & $-2.4$ & $-2.4$ & $-3.1$ \\
 $xy$   & $-2.2$ & $-2.3$ & $-2.5$ & $-3.3$ & $-2.4$ & $-2.4$ & $-3.1$ \\
 $xz$   & $-\boldsymbol{\infty}$ & $-1.4$ & $-1.7$ & $-3.3$ & $-1.3$ & $-1.6$ & $-3.1$ \\
 $y^2$  & $-2.1$ & $-2.2$ & $-2.5$ & $-3.2$ & $-2.4$ & $-2.4$ & $-3.1$ \\
 $yz$   & $-\boldsymbol{\infty}$ & $-1.4$ & $-1.7$ & $-3.3$ & $-1.3$ & $-1.6$ & $-3.1$ \\
 $z^2$  & $-1.9$ & $-2.0$ & $-2.2$ & $-3.7$ & $-2.3$ & $-2.3$ & $-3.2$ \\\hline
 $\lambda$ & $-1.5$ & $-1.6$ & $-1.9$ & $-2.7$ & $-$ & $-$ & $-$ 
    \end{tabular}
    \caption{The median value of $\log_{10}(E_\text{rel})$ over $s$ and $r$ for each weighting scheme and each observable in $\mathcal{B}$ at $P=21$ and $N=10^6$. POT is only evaluated at $r=1$. The last row reports the absolute error of predicting the Lyapunov exponent from periodic orbits, given by \autoref{eq:eabs}. Column headers are color-coded to match \autoref{fig:scaling}.}
    \label{tab:averages}
\end{table}

There are four parameters that can impact the weights. There is (i) the cardinality of the library, $P$, as well as (ii) the choice of \textit{which} $P$ measures are used. To explore this, we randomly permute the entire collection of $125$ measures (orbits or snippets) $R-1=255$ different times, labeling these libraries $\{\mathcal{L}_r\}_{r=2}^{R=256}$; we retain $\mathcal{L}_1$ as an un-permuted library, in which UPOs are ordered by increasing symbol length. This results in $R=256$ libraries total. At each value of $P$, weights are computed using the first $P$ measures of $\mathcal{L}_r$, for each $r$. The LSW and Markov weights, which are computed from chaotic data, can also depend on how that data is collected. Namely, (iii) where the chaotic trajectory is initialized and (iv) how long the chaotic trajectory is or, more specifically, how many snapshots from the chaotic trajectory are used. To explore this, we compute $256$ different chaotic trajectories, $\{ {\bf x}^{(s)}(n\Delta t) \}_{s=1}^{256}$, where $\Delta t=2$ and $n=1,\dots,10^6$. Each chaotic sample is randomly initialized at time $t=-25$. This random initial condition is then integrated to time $t=0$ and only data with $t>0$ is retained. This ensures all trajectories are randomly initialized at time $t=0$ within the chaotic set. Averages required for calculating the LSW or Markov weights are then calculated, for each $s$, using only the first $N$ samples of the chaotic trajectory. Specifically, the LSW weights are computed via \autoref{eq:rkhsw} with
\begin{align}\label{eq:taveestimate}
    \E[a_p] &\approx \frac1N \sum_{n=1}^{N} a_p({\bf x}^{(s)}(n\Delta t)).
\end{align}
Similarly, the Markov weights are computed as
\begin{align}\label{eq:markovestimate}
w_p  &= \E\left[I({\bf x},V_p)\right]\approx \frac1N \sum_{n=1}^{N} I({\bf x}^{(s)}(n\Delta t),V_p),
\end{align}
where the Voronoi volume $V_p$ is defined as the collection of points ${\bf x}\in\mathcal{M}$ closer to subset $S_p$ than any other subset,
\begin{align}
V_p = \{ {\bf x}\in \mathcal{M}\ | \ \argmin_{q\in \mathcal{L}_r} \min_{{\bf y}\in S_q} \|{\bf x}-{\bf y}\| = p\},
\end{align}
and $I$ is the previously mentioned indicator function. Hence, for Markov weights, $w_p$ is simply the fraction of time ${\bf x}(t)$ is closer to subset $S_p$ than any other subset. 
Sampling over $P$, $r$, $s$, and $N$ allows us to eliminate most of the bias and to quantify the uncertainty of $E_\text{rel}$ and $E_\text{max}$ for each of the four weighting schemes. Of these four degrees of freedom, the uniform weighting scheme will vary only with $P$ and $r$, since it does not require chaotic data. Similarly, POT weights will vary only with $P$, as they are only computed for complete libraries (i.e. $r=1$ for specific $P$). 

\begin{figure*}
    \center
    (a)\includegraphics[width=0.4\textwidth]{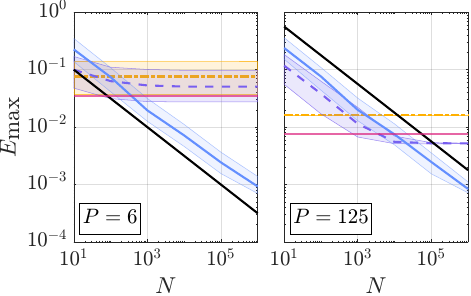}\hspace{3mm}(b)\includegraphics[width=0.4\textwidth]{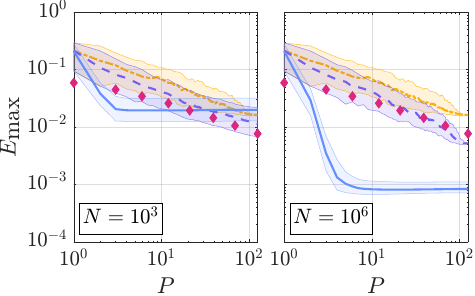}\vspace{3mm}\\
    (c)\includegraphics[width=0.4\textwidth]{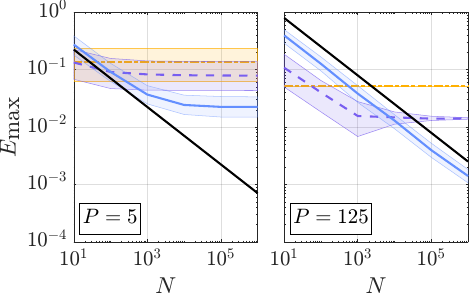}\hspace{3mm}(d)\includegraphics[width=0.4\textwidth]{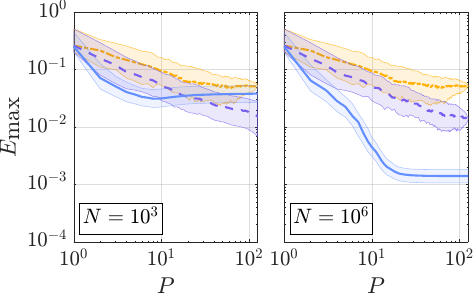}
    \caption{\label{fig:scaling}
    The median (line) and interquartile range (shaded region) of the error $E_\text{max}$  over $s$ and $r$ is plotted for orbits (panels a,b) and snippets (panels c,d). In (a,c), $P$ is fixed and $N$ is varied. In (b,d), $N$ is fixed and $P$ is varied. Line style and color correspondence is as follows: LSW (solid, blue), Markov (dashed, purple), and uniform (dot-dashed, yellow). The error of POT for complete libraries is overlaid as magenta diamonds. The scaling $E_\text{max}\sim N^{-1/2}$ is plotted as a solid black line. 
    }
\end{figure*}

\autoref{tab:averages} summarizes $E_\text{rel}$ for each method and each observable in $\mathcal{B}$ at a representative choice of $P$ and $N$. We find that, excluding errors that are identically zero, the two most accurate estimators of chaotic averages, across all observables, are the LSW estimators over UPOs and snippets. The LSW estimator over orbits almost always outperforms the estimator over snippets. 

For a number of observables, the error associated with POT is identically zero. This is because the Lorenz attractor is invariant under the transformation $(x,y,z)\mapsto(-x,-y,z)$. As a result, the infinite time average of any monomial that is odd in $x$ or $y$ is zero, e.g., $\E[x] = \E[y] = \E[xz] = \E[yz] = 0$. For every orbit $p$ of the Lorenz system, there exists a corresponding orbit $q$ such that $(x_p(t),y_p(t),z_p(t))=(-x_q(t),-y_q(t),z_q(t))$. Since these orbits must have identical stability properties, POT will assign them the same weight, $w_p=w_q$. Yet, they will produce equal but opposite time averages of any monomial that is odd in $x$ or $y$, $\E_p[a]=-\E_q[a]$. Since $p$ and $q$ must have the same symbol length, every \textit{complete} library is guaranteed to contain both orbits and their contributions must cancel, $w_p\E_p[a]+w_q\E_q[a]=0$. As a result, the POT weighting also predicts a zero infinite time average for any any monomial that is odd in $x$ or $y$. It is important to note that this property of POT is incredibly sensitive to the library being complete. If a single orbit is missed or an extraneous orbit included, this property no longer holds. 

For $a=1$, the errors associated with all non-LSW weightings are identically zero. This is because each of these methods return normalized weights, $\sum_{p=1}^P w_p = 1$, by construction; an accurate prediction of $a=1$ relies only on the weights being normalized, $1=\E[1]=\sum_{p=1}^P \E_p[1]w_p = \sum_{p=1}^P w_p$. In contrast, LSW solutions are not constrained to be normalized. Yet, LSW \textit{discovers} solutions that are normalized to $4$ significant digits, see \autoref{tab:averages}. Notably, normalization of the LSW solution can be enforced to machine precision, if desired, which we explore in \autoref{sec:uniquenessandconstraints}.

\autoref{fig:scaling} shows the dependence of $E_\text{max}$ on $N$ and $P$. 
At fixed $N$, $E_\text{max}$ decreases more rapidly in $P$ for the LSW estimate than for any other method, see \autoref{fig:scaling}(b,d); this is the case for both orbits and snippets. For $N=10^6$, the LSW estimate is often an order of magnitude more accurate than the other methods. Most importantly, LSW weights yield an accurate estimate for all considered observables from as few as $20$ orbits or snippets. 
This corresponds to accurately reproducing averages 
computed directly via \autoref{eq:taveestimate}
over $T=256\times N\Delta t = \mathcal{O}(10^9)$ time units worth of chaotic data using only $T=\mathcal{O}(10^2)$ time units worth of orbit or snippet data. 

As illustrated by \autoref{fig:scaling}(a,c), the error associated with the LSW estimate scales as $E_\text{max}\sim N^{-1/2}$ for sufficiently low $N$ and high $P$. This is most likely due to the accuracy of approximating $\E[a_p]$ via \autoref{eq:taveestimate} determined by the central limit theorem. The Markov error also scales like $E_\text{max}\sim N^{-1/2}$, and often slightly outperforms the LSW estimate, when $N<10^3$. As $N$ increases above $10^3$, the Markov estimator converges to about $1\%$ error. Hence, when chaotic data is very sparse, it seems well-motivated to use a Markov estimator of infinite-time averages. When more data is available, LSW provides noticeably more accurate estimates.

To the authors' knowledge, the convergence rates of the Markov, Uniform, and POT methods have not been directly compared in prior studies. Surprisingly, our results show that, for the range of $P$ considered here, POT does not yield meaningfully better estimates than the empirical approaches. The convergence rates, in the number of orbits/snippets, of all non-LSW orbit expansions are largely comparable. Moreover, the error of Markov weights appears to converge faster than that of POT for moderate values of $P$ considered here. When applied to snippets, Markov weights consistently outperform uniform weights; both weighting schemes appear to asymptote at large values of $P$.  

A notable advantage of POT remains that its errors do not plateau as a function of $P$; unlike the error of empirical approaches, the error of POT is proven to converge to zero as the library of orbits used becomes increasingly large and complete.
On the other hand, the LSW method is able to provide an estimate of its asymptotic accuracy, $E_\text{max}\sim N^{-1/2}$, a capability absent in both Markov and uniform weightings. For Markov weighting, it is the variance of the weights which scales like $N^{-1/2}$, rather than the accuracy of the estimate.

Finally, we consider the problem of computing the top Lyapunov exponent, $\lambda$. In particular, for the Lorenz-1963 dynamics, one finds  $\lambda = 0.90566(7)$ \cite{viswanaththesis}. It is reasonable to suppose that $\lambda$ is approximated by a weighted sum of each UPO's leading Floquet exponent, $\lambda_p$:
\begin{equation}\label{eq:lya}
    \lambda \approxeq \sum_{p=1}^P w_p \lambda_p.
\end{equation}
Indeed, when $w_p$ are POT weights and the dynamics are Axiom A, it is provable that $\sum_{p=1}^P w_p \lambda_p \to \lambda$ as $P\to\infty$ \cite{eckhardt1994}.
The final row of \autoref{tab:averages} lists the absolute error
\begin{equation}\label{eq:eabs}
E_\text{abs}(\lambda) = \left|\lambda-\sum_{p=1}^P w_p \lambda_p\right|
\end{equation}
for LSW, Markov, POT, and uniform weightings. For all $P$ examined in this paper, the LSW makes the best prediction of $\lambda$ across all methods considered, by nearly an order of magnitude. In Lorenz, the remaining Lyapunov exponents are $0$ and $-\lambda-\sigma-\beta-1$, respectively. Hence, computing both $\lambda$ and the constant observable accurately ensures that all Lyapunov exponents are computed to a similar accuracy.

\section{Uniqueness of Least Squares Weights\label{sec:uniquenessandconstraints}}

\begin{figure*}
    \center
    \subfloat[]{\includegraphics[width=0.25\textwidth]{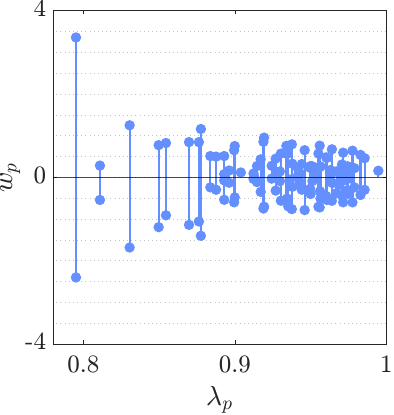}}\hspace{1mm}
    \subfloat[]{\includegraphics[width=0.255\textwidth]{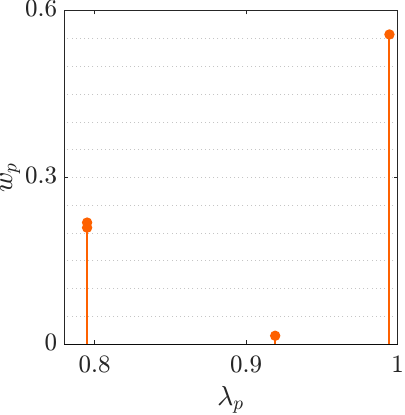}}\hspace{1mm}
    \subfloat[]{\includegraphics[width=0.265\textwidth]{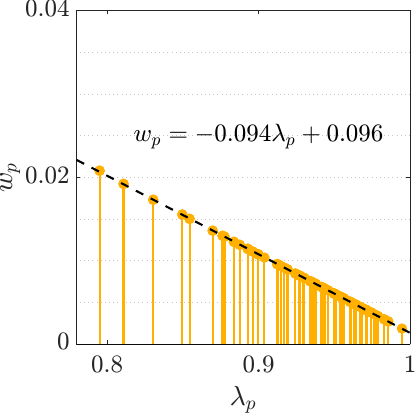}}\\
    \subfloat[]{\includegraphics[width=0.26\textwidth]{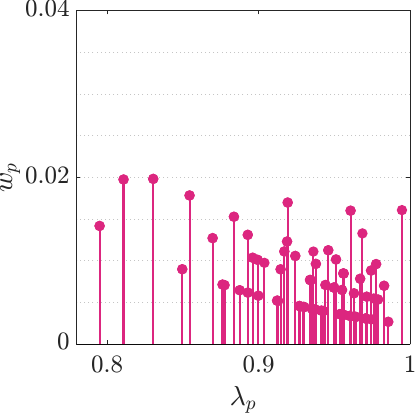}}\hspace{1mm}
    \subfloat[]{\includegraphics[width=0.26\textwidth]{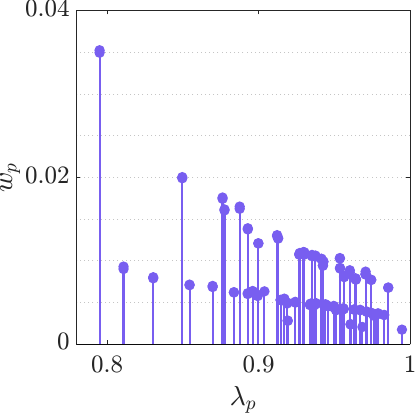}}\hspace{1mm}
    \subfloat[]{\includegraphics[width=0.272\textwidth]{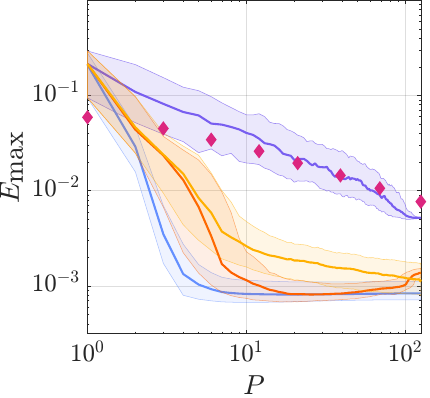}}\\
    \hspace{5mm}\subfloat[]{\includegraphics[width=0.28\textwidth]{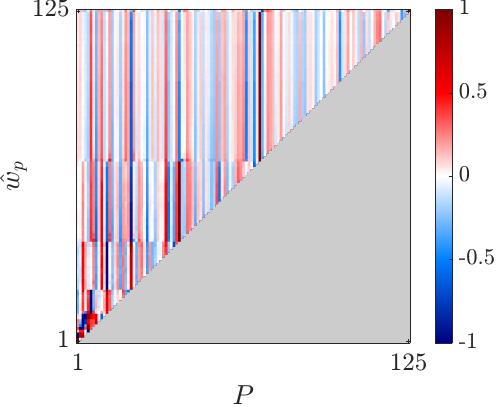}}
    \subfloat[]{\includegraphics[width=0.28\textwidth]{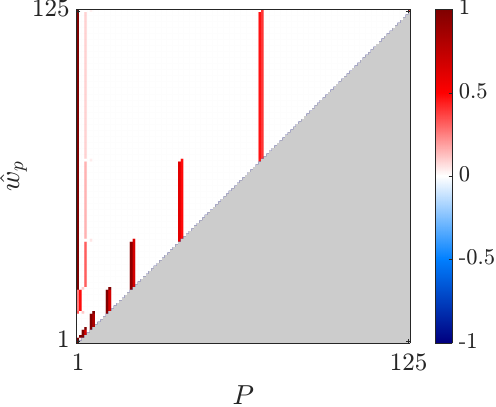}}
    \subfloat[]{\includegraphics[width=0.28\textwidth]{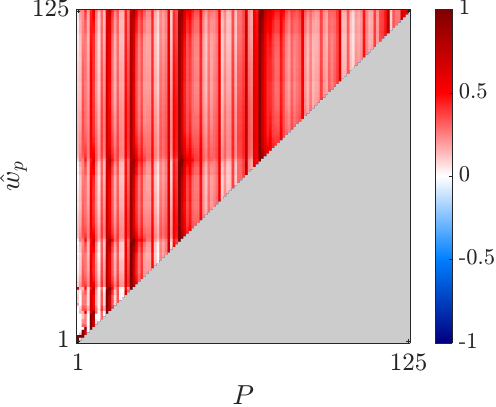}}\\
    \subfloat[]{\includegraphics[width=0.28\textwidth]{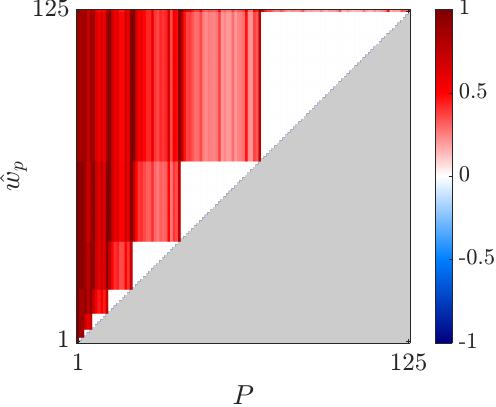}}
    \subfloat[]{\includegraphics[width=0.28\textwidth]{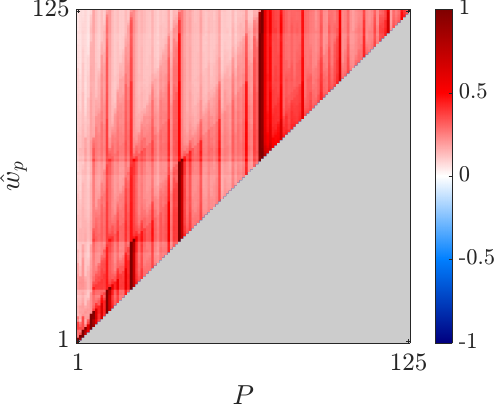}}
    \caption{\label{fig:nonuniqueness}
    A comparison of periodic orbits' weights for different methods. The results for Least Square Weighting: the dense (a) and sparse (b) solutions to \autoref{eq:sys_convex} and the solution to \autoref{eq:sys_tikhonov} (c). POT weights are shown in (d) and Markov weights in (e). In each panel (a-e), the weight $w_p$ is plotted as a function of the orbit's Floquet exponent $\lambda_p$ with $P=125$. In (c) the dashed line shows a linear fit. The error associated with each method as a function of $P$ is plotted in (f). The median error is shown as a line and the interquartile range is shaded; symbols show the POT error. The color of each bar plot matches the line color in panel (f). Panels (g-k) show the distribution of weights as a function of $P$ for each method with $r=1$ and $N=10^6$. For visualization purposes, the normalized weights $\hat{w}_p = w_{p=1}^P/\max_p|w_p|$ are shown for each $P$. Red (blue) indicates positive (negative) weight, and white indicates the weight is zero. In panel (j), POT weights are computed for the largest complete library contained within the first $P$ orbits of $\mathcal{L}_{1}$; the other orbits' weights are zeroed out. 
    }
\end{figure*}

The linear \autoref{eq:rkhsw}
has a unique solution when the matrix $A$ is non-singular, e.g., in the limit $\theta\to 0$. Convolution with a Gaussian kernel is a smoothing operation that introduces increasingly many marginal directions into $A$ as $\theta$, the smoothing parameter, is increased. In the limit $\theta \to \infty$, \autoref{eq:tutorial_continuous} limits to a single normalization constraint,
    $1 = \sum_p w_p,$
that is independent of ${\bf y}$; in this limit, the weights become severely under-constrained. Hence, \autoref{eq:rkhsw} generally admits multiple solutions. 

Marginal modes that are introduced at finite values of $\theta$ are rather interesting. Each marginal mode of $A$ approximately satisfies $\sum_p w_p = 0$ and so generates a splitting of the set of orbits into two subsets, one with positive weights $\mathcal{P} = \set{p| w_p>0}$ and another with negative weights $\mathcal{N} = \set{p| w_p<0}$. The orbits in $\mathcal{P}$ and $\mathcal{N}$ cover approximately the same regions of state space
\begin{equation}
    \sum_{p\in\mathcal{P}} w_p f_p({\bf x}) \approx \sum_{n\in\mathcal{N}} w_n f_n({\bf x}).
\end{equation}
When either $\mathcal{N}$ or $\mathcal{P}$ have only one member, this sole orbit must be a long orbit that shadows all short orbits in the other set. Hence, marginal modes can shed light on how orbits shadow each other within a chaotic system. 


When weights are both positive ($w_p\geq 0$, for all $p$) and normalized ($\sum_p w_p=1$), each $w_p$ is interpretable as the probability that a randomly sampled chaotic state is well-approximated by a state from orbit $p$. This makes positivity and normalization enticing properties for weights to satisfy. Normalization is required of any collection of weights that accurately predict a constant observable. It is less clear whether positivity is required of the weights. When the sets $\set{S_p}$ are non-overlapping, it follows that any weights satisfying \autoref{eq:genexpansion} must be positive. The SRB measure is a probability measure, meaning that $\E[I({\bf x},A)] \geq 0$ for any set $A\subset \mathcal{M}$. Applying \autoref{eq:genexpansion} to the indicator function $I({\bf x},S_q)$, it follows that 
$$\sum_{p=1}^P w_p \E_p[I({\bf x},S_q)] = w_q = \E[I({\bf x},S_q)] \geq 0 \qquad q=1,\dots,P.$$
However, if the sets $\set{S_p}$ are overlapping, then the optimal weights can be negative, even when both $\mu$ and all $\mu_p$ are probability measures. As a simple example, consider expanding the uniform probability measure on a set of four points $\mu = \set{1/4,1/4,1/4,1/4}$ in terms of probability measures $\mu_1 = \set{1/3,1/3,1/3,0}$, $\mu_2 = \set{0,1/3,1/3,1/3}$, and $\mu_3 = \set{0,1/2,1/2,0}$. It follows that 
$$\mu = \frac{3}{4}\mu_1+\frac{3}{4}\mu_2-\frac{1}{2}\mu_3,$$
where the weight associated with $\mu_3$ is negative.

Markov weights are normalized and positive by definition. POT weights are guaranteed to be normalized, however, it is unproven whether or not POT weights are positive. Empirically, we find that POT weights are positive when computed on complete libraries, but can produce negative weights when evaluated on incomplete libraries. The Tikhonov-regularized solution used in \autoref{sec:results} does not enforce positivity or normalization and is not guaranteed to produce weights that can be interpreted as probabilities. As mentioned, we find that the solution found by Tikhonov regularization happens to be normalized to nearly four significant digits. Moreover, we observe that the Tikhonov-regularized solution becomes increasingly positive (and remains approximately normalized) as the regularization parameter $\alpha$ is increased. However, the residual of the original least-squares problem, $\|A\vec{w}-\vec{b}\|$, increases with increasing $\alpha$. The weights resulting from Tikhonov regularization with $\alpha=10^{-10}$ and $P=125$ are plotted as a function of orbit stability in \autoref{fig:nonuniqueness}(a). The distribution of $w_p$ is plotted as a function of $P$ over the library $\mathcal{L}_1$ in \autoref{fig:nonuniqueness}(g). The distribution of $w_p$ is largely smoothly varying as a function of $P$, with notable discontinuities at values of $P$ corresponding to complete libraries.  

Is it possible to constrain the LSW solution to be both positive and normalized. That is, we may seek minima of the function
\begin{align}\label{eq:sys_convex}
&V=\|A\vec{w} - \vec{b}\|_2^2\\
\text{subject to}\qquad &w_p \geq 0,\hspace{3mm} \sum_{p=1}^P w_p = 1.\nonumber
\end{align}
We do this two ways. First, we approximate a solution to \autoref{eq:sys_convex} using MATLAB's built-in \verb|lsqnonneg| function and normalizing the output. The output of \verb|lsqnonneg| happens to be normalized about four decimals places for all $P$, and the normalization we apply afterwards does not greatly perturb the solution away from a minimum of $V$. More importantly, as a consequence of MATLAB's implementation of \verb|lsqnonneg|, these solutions are sparse.
We found this to be the most reliable way to find a sparse solution. The resulting weights are plotted in \autoref{fig:nonuniqueness}(b,h). In terms of the standard labeling convention for symbolic dynamics \cite{viswanath2003}, the sparse solution prefers to retain orbit $AB$, the shortest orbit of the Lorenz dynamics, and highly asymmetric orbits of the form $A^nB$ and $B^n A$ for integer $n$. 

We also use MATLAB's built-in \verb|fmincon| function to minimize $V$ directly. The solution found by \verb|fmincon| does depend on initial conditions. In particular, the weights resulting from minimization using the initial condition $w_p=1/P$ are plotted in \autoref{fig:nonuniqueness}(c,i). 

It has been posited before that $w_p$ should be inversely proportional to an orbit's stability \cite{kazantsev1998, chandler2013}; intuitively, the more unstable an orbit is compared to others, the less frequently it should be visited by the chaotic trajectory, and the shorter the duration of each visit. Interestingly, almost all weighting schemes (save the sparse one) indicate that $|w_p|$ tends to decrease with an orbit's Floquet exponent $\lambda_p$, although the dependence is generally nonmonotonic. Surprisingly, the solution found by \verb|fmincon| yields a simple linear relationship between $w_p$ and $\lambda_p$ (rather than $\lambda_p^{-1}$, as has been suggested previously). Note that such linear relationship is only found for complete libraries (e.g., $P=125$).  

Perhaps the most surprising result of this study is that \autoref{eq:genexpansion} admits multiple, qualitatively distinct distributions of weights that yield accurate estimates of time-averages, even when $P$ is small (see \autoref{fig:nonuniqueness}(a-c,f)). To the authors' knowledge, quantitative evidence of this has not been presented in the literature thus far. Most impressively, the sparse solution found by \verb|lsqnonneg| achieves comparable accuracy to the unconstrained LSW solution and returns weights that may be interpreted as probabilities, while never requiring more than $4$ orbits in the expansion. 

\section{Conclusion}

To summarize, we have derived a general method for performing least squares regression on linear relations of measures, and used it to construct an accurate approximation of infinite-time averages for dynamical systems in terms of a weighted average over as few as 4 periodic orbits or 10 snippets of chaotic trajectory. Thus, the approach introduced here works almost equally well whether or not any periodic orbits are known. 

The weights are given by a least squares solution to a linear system of equations, and, unlike those predicted by period orbit theory, are quite straightforward to compute. As demonstrated here, enforcing constraints on the weights, such as positivity and normalization, is also simple using off-the-shelf linear programming algorithms.

We have also identified the subspace of observables over which the Least Squares Weighting yields exact predictions. Even outside of this subspace, the accuracy of the method tends to be an order of magnitude better than state-of-the-art methods, and more quickly convergent in the number of orbits or snippets used. Surprisingly, the method discovered \textit{multiple} weightings that outperform state-of-the-art methods. 


Given that the Least Squares Weighting is based on Kriging, it is expected to scale quite well to high- and even infinite-dimensional state spaces, such as those needed to describe fluid, plasma, or optical turbulence, although this has yet to be verified. In such applications, few if any periodic orbits are typically known, so the ability of our approach to work when applied to few orbits and snippets of chaotic trajectory makes it uniquely appealing. 

\begin{acknowledgments}
We gratefully acknowledge financial support from the National Science Foundation under {Grant~No.~2032657.}
\end{acknowledgments}

\section*{Data Availability Statement}

The code and data that support the findings of
this study are openly available on Github \cite{repository}. DOI will be added upon publication. 
The algorithm used for computing POT weights is described within the article's supplementary material.

\appendix

\section{Computing POT Weights Numerically}

The theoretical foundations of POT are summarized by \citet{cvitanovic1999}, and Chaosbook provides a detailed reference text on the subject \cite{chaosbook} .
\citet{eckhardt1994} applied POT to the Lorenz 1963 system. Here, we provide a short summary of the essential steps. The following algorithm works for systems with symbolic dynamics, however, averages can still be computed from POT without one \cite[Ch. 23.7]{chaosbook}. \citet{viswanath2003} provides a good (and visual) description of symbolic dynamics in the context of Lorenz. 

Given a library of periodic orbits which contains all orbits up to symbolic length $n$, the temporal average is approximated as 
\begin{align}\label{app:a}
    \mathbb{E}[a] \approx  -\frac{\partial_\beta F_n(s,\beta)}{\partial_s F_n(s,\beta)}\Bigg|_{\substack{\ s=s_0,\\ \beta=0\ }},
\end{align}
where $F_n$ approximates a spectral determinant\cite[Ch. 22]{chaosbook} and $s_0$ is the real zero of $F_n$, i.e., $F_n(s_0,0)=0$. Equality is reached in \autoref{app:a} in the limit $n\to \infty$. While POT does not define the weights $w_p$ explicitly, they may be computed from \autoref{app:a}:
\begin{align}\label{app:w}
    w_p =  -\frac{\partial_{\mu_{p}}\partial_\beta F_n(s,\beta)}{\partial_s F_n(s,\beta)}\Bigg|_{\substack{\ s=s_0,\\ \beta=0\ }}.
\end{align}

The spectral determinant can be computed using the trace coefficients \cite[Ch. 23.2.2]{chaosbook}
\begin{align}
C_j = - \sum_{n_pr = j} \frac{1}{r}\frac{e^{-rT_p \left(s-\beta\mu_{p}(a)\right) }}{\lvert\det(1-(M_p)^r)\rvert},
\end{align}
where the sum is over all prime orbits $p$ such that $j=rn_p$, where $r=1,2,\dots$ and $n_p$ is the symbol length of orbit $p$. Notice that $C_j$ sums over all periodic orbits of length $j$ \textit{and} repeats of shorter orbits which, when repeated $r$ times, have symbol length $j$. $M_p$ is the transverse monodromy matrix \cite[Ch. 21.2.2]{chaosbook} of orbit $p$, over one period. For the Lorenz system, $\det(1-(M_p)^r) \approx  1-e^{rT_p\lambda_p}$, where $\lambda_p$ and $T_p$ are the top Floquet exponent and period of orbit $p$, respectively.

The spectral determinant is then
\begin{align}\label{eq:qSum}
F_n(s,\beta) &= 1 + \sum_{j=1}^n Q_j
\end{align}
with the help of a recurrence relation
\begin{align}\label{eq:qRecur}
Q_j & = C_j + \sum_{i = 1}^{j-1} \frac{j-i}{j}C_{j-i}Q_i,
\end{align}
where $ Q_1= C_1$.
The partial derivatives of $F_n$ are computed via the chain rule with
\begin{align}
\partial_s C_j\Big|_{\beta=0 } &=  {\sum_{n_pr = j}} \frac{e^{-rT_ps}}{\lvert 1-e^{rT_p\lambda_p}\rvert}T_p,\label{eq:chainS}\\
\partial_\beta C_j\Big|_{\beta=0 } &= - {\sum_{n_pr = j}} \frac{e^{-rT_ps}}{\lvert 1-e^{rT_p\lambda_p}\rvert}T_p\langle a\rangle_p, \label{eq:chainB}\\
\partial_{\mu_p}\partial_\beta C_j \Big|_{\beta=0 }&= - \sum_{n_pr = j} \frac{e^{-rT_ps}}{\lvert 1-e^{rT_p\lambda_p}\rvert}T_p,\label{eq:chainW}\\
\partial_{\mu_p} C_j \Big|_{\beta=0 }&= 0.
\end{align}
In \autoref{eq:chainW}, note that $p$ is fixed and the sum runs only over $r$. 
The root $s_0$ is found using Newton's method by iterating the relation
\begin{align}
        s_0 = s_0 - \frac{F_n(s_0,0)}{\partial_s F_n(s_0,0)}
\end{align} 
until $|F_n(s_0,0)| < \epsilon$ using $s_0=0$ as the initial guess. We use $\epsilon= 10^{-8}$. 

\bibliography{references}

\end{document}